\documentclass[a4paper,10pt]{article}

\usepackage{csquotes}

\usepackage{amsmath}
\usepackage{amsfonts}
\usepackage{amsthm}
\usepackage{enumerate}
\usepackage{color}
\usepackage{fancyref}
\usepackage{color}
\usepackage{verbatim}
\usepackage{amsbsy}
\usepackage{multirow}
\usepackage{mathrsfs}
\usepackage{enumitem}
\usepackage{mathtools}
\usepackage{epsfig}
\usepackage[strict]{changepage}
\usepackage{amssymb}
\usepackage[UKenglish]{babel}
\usepackage{cite}

\numberwithin{equation}{section}

\theoremstyle{plain}
\newtheorem{theorem}{Theorem}[section]
\newtheorem{lemma}[theorem]{Lemma}
\newtheorem{proposition}[theorem]{Proposition}

\theoremstyle{remark}
\newtheorem{remark}[theorem]{Remark}
\newtheorem{example}[theorem]{Example}
\theoremstyle{definition}

\DeclareMathOperator\real{Re}
\DeclareMathOperator\imag{Im}
\renewcommand\Re{\real}
\renewcommand\Im{\imag}
\newcommand\rd{\mathrm{d}}
\newcommand\CC{\mathbb{C}}
\newcommand\NN{\mathbb{N}}
\newcommand\RR{\mathbb{R}}

\newcommand\mr{\mathring}
\newcommand\wh{\widehat}
\newcommand\wt{\widetilde}
\DeclareMathOperator\spn{span}

\newcounter{counter_a}
\newenvironment{myenum}{\begin{list}{{\rm(\roman{counter_a})}}%
{\usecounter{counter_a}
\setlength{\itemsep}{0.5ex}\setlength{\topsep}{0.7ex}
\setlength{\leftmargin}{5ex}\setlength{\labelwidth}{5ex}}}{\end{list}}

\newcounter{counter_b}
\newenvironment{myenuma}{\begin{list}{\textup{(\alph{counter_b})}}%
{\usecounter{counter_b}
\setlength{\itemsep}{0.5ex}\setlength{\topsep}{0.7ex}
\setlength{\leftmargin}{5ex}\setlength{\labelwidth}{5ex}}}{\end{list}}

\newcounter{counter_c}
\newenvironment{myenum_A}{\begin{list}{\textup{(A\arabic{counter_c})}}%
{\usecounter{counter_c}
\setlength{\itemsep}{0.5ex}\setlength{\topsep}{0.7ex}
\setlength{\leftmargin}{7ex}\setlength{\labelwidth}{7ex}}}{\end{list}}

\newcounter{counter_d}
\newenvironment{myenum_EF}{\begin{list}{\textup{(EF\arabic{counter_d})}}%
{\usecounter{counter_d}
\setlength{\itemsep}{0.5ex}\setlength{\topsep}{0.7ex}
\setlength{\leftmargin}{8ex}\setlength{\labelwidth}{8ex}}}{\end{list}}

\newenvironment{myenum_P}{\begin{list}{\textup{(P\arabic{counter_d})}}%
{\usecounter{counter_d}
\setlength{\itemsep}{0.5ex}\setlength{\topsep}{0.7ex}
\setlength{\leftmargin}{7ex}\setlength{\labelwidth}{7ex}}}{\end{list}}

\begin{document}

\title{Discrete Fragmentation Equations with Time-Dependent Coefficients}
\date{}
\author{Lyndsay Kerr, Wilson Lamb and Matthias Langer}

\maketitle

\begin{center}
\textit{Dedicated to Jerry Goldstein on the occasion of his 80$^{\textit{th}}$ birthday}
\end{center}

\medskip

\begin{abstract}
\noindent
We examine an infinite, linear system of ordinary differential equations
that models the evolution of fragmenting clusters, where each cluster is assumed
to be composed of identical units.  In contrast to previous investigations into
such discrete-size fragmentation models,  we allow the fragmentation coefficients
to vary with time.
By formulating the initial-value problem for the system as a non-autonomous
abstract Cauchy problem, posed in an appropriately weighted $\ell^1$ space,
and then applying results from the theory of evolution families,
we prove the existence and uniqueness of physically relevant, classical solutions
for suitably constrained coefficients.
\\[1ex]
\\[1ex]
\end{abstract}

\section{Introduction}

Fragmentation is a commonly observed phenomenon in various physical processes
such as polymer degradation, liquid droplet breakup, and the crushing and grinding of rocks.
Deterministic models of fragmentation are usually based on the simplifying assumption
that the fragmenting objects can be distinguished by means of a single `size' variable,
such as mass. When this size variable is permitted to take any positive value,
the resulting model describing the continuous time evolution of the system of
fragmenting objects typically takes the form of a linear integro-differential equation
which is referred to as the continuous (size) fragmentation equation;
see \cite[\S 2.2.2]{banasiak2019analytic} and
also the discussion in Section~\ref{sect:concluding_remarks}.
However, when describing the fragmentation of clusters that are assumed to be
comprised of a finite number of identical fundamental particles, it is clearly
more appropriate to use a discrete-size variable.
It is the discrete-size case that we examine in this paper, and, for convenience,
we adopt the polymer-based terminology that is frequently used when dealing
with discrete-size fragmentation.  Consequently, the fundamental particle is
referred to as a monomer, and an $n$-mer is then a cluster consisting of $n$ monomers.
By suitable scaling, a monomer can be assumed to have unit mass in which case
an $n$-mer has mass~$n$.

In terms of the number density, $u_n(t)$, of $n$-mers at time $t$, the evolution
of the system of fragmenting clusters is described by the infinite-system of
linear ordinary differential equations (ODEs)
\begin{equation}\label{NA frag system}
\begin{split}
  u_n'(t) &= -a_n(t)u_n(t)+\sum_{j=n+1}^{\infty} a_j(t)b_{n,j}(t)u_j(t),
  \quad t\in(0,T],\ n \in \NN; \\
  u_n(0) &= \mathring{u}_n, \quad n \in \NN,
\end{split}
\end{equation}
where $T>0$, and $\mr{u}$ is an initial density sequence.
At each time $t \in [0,T]$, the coefficients $a_n(t)$ and $b_{n,j}(t)$ represent,
respectively, the rate at which $n$-mers are lost due to fragmentation (when $n \ge 2$),
and the average number of $n$-mers that are produced when a larger $j$-mer fragments.
Throughout, we assume that the fragmentation coefficients satisfy the following
natural physical constraints:
\begin{myenum_A}
\item
  $a_n(t) \ge 0, \; \forall t \in [0,T]$ and $n\in\NN$,
\item
  $b_{n,j}(t) \ge 0, \;  \forall t \in [0,T]$ and $n,j \in\NN$,
with $b_{n,j}(t) = 0$ if $j \le n$.
\end{myenum_A}
As monomers cannot fragment to produce smaller clusters, the case $a_1(t) > 0$
signifies a depletion in the number of monomers due to some other mechanism;
see \cite{caiedwardshan1991} and \cite{smith2012discrete}.

The total mass of all clusters at time $t \in [0,T]$ is given by the first moment,
$M_1(u(t))$, of the sequence of densities $u(t)=(u_n(t))_{n=1}^\infty$, where
\begin{equation}\label{total mass}
  M_1\bigl(u(t)\bigr) \coloneqq \sum_{n=1}^{\infty} nu_n(t).
\end{equation}
On representing the total mass of daughter clusters produced from the
fragmentation of a $j$-mer at time $t$ by
\begin{equation}\label{local mass conservation lambda}
  \sum_{n=1}^{j-1} nb_{n,j}(t) = \bigl(1-\lambda_j(t)\bigr)j, \qquad j=2,3,\ldots,
\end{equation}
where each $\lambda_j$ is a real-valued function,
a formal calculation establishes that if $u(t) = (u_n(t))_{n=1}^\infty$
is a solution of \eqref{NA frag system}, then
\begin{equation}\label{massode}
  \frac{\rd}{\rd t}\Bigl(M_1\bigl(u(t)\bigr)\Bigr)
  = - a_1(t)u_1(t) - \sum_{j=2}^\infty j \lambda_j(t) a_j(t)u_j(t), \quad t\in(0,T].
\end{equation}
The expression in \eqref{massode} gives the rate at which mass may be lost or gained
from the system of fragmenting clusters, and also shows, at least formally,
that there is no change in the total mass during the fragmentation process
when $a_1(t)=0$ and $\lambda_j(t)=0$ for all $j=2,3,\ldots$, and $t \in [0,T]$.
Let us note that we do not impose any sign restrictions on $\lambda_j(t)$, i.e.\
we also allow mass to be gained in a fragmentation event,
in which case $\lambda_j(t)<0$.

In contrast to the continuous fragmentation equation, where several investigations,
such as \cite{mclaughlin1997NA} and \cite{melzak1957scalar}, have dealt with
time-dependent coefficients, previous investigations into \eqref{NA frag system}
appear to have considered only the case when all fragmentation coefficients are
time-independent.
The approach used in \cite{banasiak2012global, mcbride2010strongly,smith2012discrete}
to analyse the constant-coefficient fragmentation system is to formulate the
initial-value problem as an autonomous abstract Cauchy problem (ACP), posed in
an appropriate Banach lattice.
Conditions on the coefficients are then determined under which the ACP has a
unique classical solution that can be expressed in terms of a positive $C_0$-semigroup
of contractions (i.e.\ a substochastic semigroup), usually referred to as the
fragmentation semigroup.

In a recent paper \cite{kerr2020discrete}, we also applied this semigroup-based strategy,
working within the framework of weighted $\ell^1$ spaces of the form
\begin{equation}\label{weighted l^1 space}
  \ell_w^1 \coloneqq \biggl\{f=(f_n)_{n=1}^{\infty}: f_n \in\RR,\, \forall n \in\NN,
  \text{ and }  \|f\|_w \coloneqq \sum_{n=1}^{\infty} w_n|f_n|<\infty\biggr\},
\end{equation}
where $w_n>0$, $n=1,2,\ldots$.
Each space $\ell_w^1$ is a real Banach lattice with positive cone
\begin{equation}\label{poscone}
  (\ell_w^1)_+ \coloneqq \bigl\{f=(f_n)_{n=1}^\infty \in \ell_w^1 : f_n \ge 0, \,
  \forall n \in\NN\bigr\}.
\end{equation}
Moreover, since $\|\cdot\|_w$ is additive on $(\ell_w^1)_+$, $\ell_w^1$ is an $AL$-space;
see \cite[Definition~2.56]{banasiak2006perturbations}.
By allowing general weights, $w=(w_n)_{n=1}^{\infty}$ with $w_n > 0$ for all $n \in\NN$,
and not just the specific cases of $w_n = n$ and $w_n = n^p$, $p > 1$, used, respectively,
in \cite{mcbride2010strongly,smith2012discrete} and \cite{banasiak2012global},
we were able to establish the following results.

Firstly, given \emph{any} fragmentation coefficients satisfying the time-independent versions
of (A1) and (A2), it transpires that it is \emph{always} possible to determine
a weight which will guarantee the existence of a substochastic fragmentation semigroup
on $\ell_w^1$.
More precisely, from \cite[Theorem~3.4]{kerr2020discrete}, there exists a
substochastic fragmentation semigroup on $\ell_w^1$ whenever $w_n \ge n$
for all $n\in\NN$, and there exists $\kappa \in (0,1]$ such that
\begin{equation}\label{weight assumption for A frag}
  \sum_{n=1}^{j-1} w_nb_{n,j} \le \kappa w_j,\; \forall j =2,3,\ldots.
\end{equation}

Secondly, under the more restrictive condition that $0<\kappa<1$,
it is shown \cite[Theorem~5.2]{kerr2020discrete} that the fragmentation semigroup
is analytic when defined on the complexification of $\ell_w^1$.
Details on the process of complexification of a real Banach lattice can be found in
\cite[\S 2.2.5]{banasiak2006perturbations}.
For the particular case of $\ell_w^1$, the process leads simply to the
complex Banach lattice defined as in \eqref{weighted l^1 space} but now for
complex sequences $f = (f_n)_{n=1}^\infty$, $f_n\in\CC$.
The partial order in this complex Banach lattice is given by
\[
  f = (f_n)_{n=1}^\infty \le g = (g_n)_{n=1}^\infty \iff
  \Re f_n \le \Re g_n \text{ and } \Im f_n = \Im g_n, \ \forall n \in \NN,
\]
and this ensures that the corresponding positive cone is also given by \eqref{poscone}.
Clearly, for any given coefficients $b_{n,j}$, a sequence $(w_n)_{n=1}^\infty$
can be constructed iteratively such that $w_n \ge n$
and \eqref{weight assumption for A frag} is satisfied for some $\kappa \in (0,1)$;
see \cite[Theorem~5.5]{kerr2020discrete}.

Our aim in the current paper is to exploit the above result on the analyticity
of the fragmentation semigroup in weighted $\ell^1$ spaces to determine sufficient conditions
under which the non-autonomous fragmentation system \eqref{NA frag system} is well posed.
In keeping with the semigroup approach used for the autonomous system, the strategy
we adopt involves the application of the theory of evolution families,
an account of which can be found in the seminal books on semigroups of operators
by Goldstein \cite{goldstein2017semigroups} and Pazy \cite{pazy1983semigroups}.
Such families have been employed in the analysis of a variety of linear,
non-autonomous evolution equations, such as the time-dependent coefficient versions
of the continuous integro-differential fragmentation equation \cite{mclaughlin1997NA}
and the Black--Scholes equation \cite{ruiz2020chaotic}.

In Section~\ref{sect:prelim}, we give some prerequisite information on non-autonomous ACPs
and strongly continuous evolution families, and then express \eqref{NA frag system}
in the form of a non-autonomous ACP posed in an $\ell_w^1$ space.
In Section~\ref{sect:well-posedness}, this abstract formulation of \eqref{NA frag system}
is shown to be well posed when the weight $w$ and the fragmentation coefficients
are suitably constrained.
In Section~\ref{sect:asymptotics} we consider the asymptotic behaviour of solutions
as $t\to\infty$.  In particular, under the assumption of mass conservation
we prove that solutions converge to a monomeric state with an explicit exponential rate.
Finally, in Section~\ref{sect:concluding_remarks}, some potential extensions to the
work presented here are discussed.

\section{Preliminaries and Abstract Formulation}
\label{sect:prelim}

To enable an approach based on the theory of evolution families to be applied
to the fragmentation system, the initial-value problem \eqref{NA frag system}
must first be recast as a non-autonomous ACP.
It turns out to be useful to allow arbitrary initial times.
So, for fixed $T>0$ we consider the family of initial-value problems
\begin{equation}\label{eq general NA ACP}
  u'(t)=G(t)u(t), \quad  t\in(s,T]; \qquad u(s)=\mr{u},
\end{equation}
where $s\in[0,T)$.
Moreover, for each $t\in[0,T]$, $G(t)$ is a linear operator that
maps $D(G(t)) \subseteq X$ into $X$, where $X$ is a Banach space,
and $\mr u\in X$ is an initial value.
The aim is to determine conditions on $G(t)$ which ensure
that \eqref{eq general NA ACP} has a unique solution $u:[s,T] \to X$
that can be expressed in terms of a (strongly continuous) evolution family
(or evolution system), which, from \cite[Definition~5.5.3]{pazy1983semigroups},
is a two-parameter family of bounded linear operators $(U(t,s))_{0\le s \le t \le T}$,
on $X$ satisfying
\begin{myenum_EF}
\item
  $U(s,s) = I,\;
  U(t,r)U(r,s) = U(t,s) \;\; \text{for }\; 0\le s \le r \le t \le T$,
\item
  $(t,s) \mapsto U(t,s)$ is strongly continuous for $0 \le s \le t \le T$.
\end{myenum_EF}

A function $u$ is a \emph{classical solution} of \eqref{eq general NA ACP}
if $u \in C([s,T],X) \cap C^1((s,T],X)$, $u(t) \in D(G(t))$ for all $t \in(s,T]$,
and \eqref{eq general NA ACP} is satisfied.
The following result, which is a slightly modified version
of \cite[Theorem~5.6.8]{pazy1983semigroups}
(see also \cite[Theorem~5.6.1]{pazy1983semigroups}), gives sufficient conditions
on the operators $G(t)$, $0 \le t \le T$, for the existence of a unique classical solution
to \eqref{eq general NA ACP}, and also highlights the key role played by
evolution families and analytic semigroups.

\begin{theorem}\label{thm pazy NA}
Let $X$ be a complex Banach space.  For each $t \in [0,T]$, where $T>0$,
let $G(t)$ be the generator of an analytic semigroup, $(S_t(\tau))_{\tau\ge0}$, on $X$.
Assume that the following conditions are satisfied.
\begin{myenum_P}
\item
  The domain $D(G(t))\eqqcolon\mathcal{D}$ of $G(t)$ is independent of $t \in [0,T]$.
\item
  For $t \in [0,T]$, the resolvent $R(\lambda,G(t))\coloneqq(\lambda I-G(t))^{-1}$
  exists for all $\lambda\in\CC$ with $\Re\lambda\ge0$ and there is a constant $M$
  such that
  \[
    \big\|R(\lambda,G(t))\big\| \le \frac{M}{|\lambda|+1} \qquad
    \text{for all} \ t \in [0,T], \text{ and } \lambda \in \CC: \Re\lambda\ge0.
  \]
\item
  There exist constants $L$ and $\sigma \in (0,1]$ such that
  \[
    \big\|\bigl(G(t)-G(s)\bigr)G(\tau)^{-1}\big\| \le L|t-s|^\sigma \qquad
    \text{for} \ s,t,\tau\in[0,T].
  \]
\end{myenum_P}
Then there exists an evolution family $(U(t,s))_{0\le s\le t\le T}$ such that,
for every $s\in[0,T)$ and $\mr{u} \in X$, the non-autonomous ACP
\eqref{eq general NA ACP} has a unique classical solution which is given
by $u(t)=U(t,s)\mr{u}$.
\end{theorem}

Note that assumption (P2) implies that the semigroup $(S_t(\tau))_{\tau\ge0}$
is analytic.  We included it explicitly in the formulation to emphasise the
importance of analyticity.
Further, note that $U(t,s)$ maps $X$ into $\mathcal{D}$ if $s<t$ because $u(t)=U(t,s)\mr{u}$
is a classical solution of \eqref{eq general NA ACP}.

The abstract formulation of the fragmentation system \eqref{NA frag system}
is posed in the complex Banach lattice version of the space $\ell_w^1$
from \eqref{weighted l^1 space},
which will also be denoted by $\ell_w^1$.
The weight $w = (w_n)_{n=1}^\infty$ is assumed to satisfy
\begin{myenum_A}
\setcounter{counter_c}{2}
\item
  $w_n \ge n$ for all $n\in\NN$,
\item
  $\exists \kappa \in (0,1): \sum\limits_{n=1}^{j-1} w_nb_{n,j}(t) \le \kappa w_j$
  for all $t\in[0,T]$ and $j=2,3,\ldots$.
\end{myenum_A}
It follows from (A3) that $\ell_w^1$ is continuously embedded in the space
\[
  X_{[1]} \coloneqq\biggl\{f=(f_n)_{n=1}^{\infty}: f_n \in \RR,\,  \forall n \in \NN,
  \text{ and }  \|f\|_{[1]} \coloneqq\sum\limits_{n=1}^{\infty} n|f_n|<\infty\biggr\},
\]
which is often referred to as the first moment space since
\[
  \|f\|_{[1]} = M_1(f) \qquad\text{for all} \ f \in (X_{[1]})_+.
\]
Moreover, on defining the bounded linear functional $\phi_w$ on $\ell_w^1$ by
\begin{equation}\label{def_phi_w}
  \phi_w\bigl((f_n)_{n=1}^\infty\bigr) \coloneqq \sum_{n=1}^\infty w_nf_n,
  \qquad (f_n)_{n=1}^\infty \in \ell_w^1,
\end{equation}
it is clear that
\[
  \phi_w(f)=\|f\|_w \qquad \text{for all} \ f\in\bigl(\ell_w^1\bigr)_+,
\]
and so $\phi_w$ coincides with the norm $\|\cdot \|_w$ on the positive cone.

\begin{remark}
\rule{0ex}{1ex}
\begin{myenum}
\item
If $b_{n,j}$ is bounded on $[0,T]$ for all $n,j\in\NN$, $n<j$, then one can
construct a sequence $(w_n)_{n=1}^\infty$ iteratively such that (A3) and (A4)
are satisfied.
\item
It can be shown in a similar way as in \cite[Theorem~5.5 and Lemma~5.4]{kerr2020discrete}
that, when
\begin{equation}\label{no_mass_gain}
  \sum_{n=1}^{j-1} nb_{n,j}(t) \le j, \qquad j=2,3,\ldots,
\end{equation}
(or equivalently $\lambda_j(t)\ge0$ in \eqref{local mass conservation lambda}),
then one can choose the sequence $(w_n)_{n=1}^\infty$ such that it grows at most
exponentially.  The condition \eqref{no_mass_gain} means that
the total mass does not grow in each fragmentation event.
\end{myenum}
\end{remark}

\medskip

\begin{example}\label{example_powers}
Consider the case when
\[
  b_{n,j}(t) \equiv b_{n,j} = \beta_n\, \zeta_j, \qquad n,j\in\NN,\; n<j,\; t\in[0,T],
\]
where $\beta_n = n^\nu$ with $\nu\ge-1$. Under the assumption that mass is
conserved during each fragmentation event
(i.e.\ $\sum_{n=1}^{j-1} nb_{n,j} = j$ for $j=2,3,\ldots$), we then obtain
\[
  \zeta_j = \frac{j}{\sum\limits_{l=1}^{j-1}l^{\nu+1}}\,.
\]
In this case, we can show that (A4) is satisfied with $w_n=n^p$ for some $p\ge1$.
To this end, we consider
\[
  \frac{1}{w_j}\sum_{n=1}^{j-1}w_n b_{n,j}
  = \frac{1}{j^p}\cdot\frac{j}{\sum\limits_{l=1}^{j-1}l^{\nu+1}}
  \sum_{n=1}^{j-1}n^{p+\nu}
\]
for $j=2,3,\ldots$.
We can estimate the two sums with integrals:
\begin{align*}
  \sum_{n=1}^{j-1}n^{p+\nu} &\le \int_1^j x^{p+\nu}\,\rd x
  \le \frac{1}{p+\nu+1}j^{p+\nu+1},
  \\[1ex]
  \sum_{l=1}^{j-1}l^{\nu+1} &\ge \int_0^{j-1} x^{\nu+1}\,\rd x
  = \frac{1}{\nu+2}(j-1)^{\nu+2}
  \ge \frac{1}{\nu+2}\Bigl(\frac{j}{2}\Bigr)^{\nu+2},
\end{align*}
which yields
\begin{equation}\label{est_example}
  \frac{1}{w_j}\sum_{n=1}^{j-1}w_n b_{n,j}
  \le \frac{(\nu+2)2^{\nu+2}}{p+\nu+1}\,.
\end{equation}
For fixed $\nu\ge-1$ we can choose $p\ge1$ such that the
right-hand side of \eqref{est_example} is strictly less than 1,
which shows that (A4) is satisfied.
\hfill$\lozenge$
\end{example}

\begin{remark}
It is worth noting that analogous separable coefficients which take the
form $b(x,y) = \beta(x)\,\zeta(y)$, $0 < x < y$, have been considered in
investigations into the continuous, mass-conserving, autonomous fragmentation equation.
In particular, the case $\beta(x) = x^\nu$ and
\[
  \zeta(y) = (\nu + 2)y^{-\nu - 1} = \frac{y}{\int_0^y y^{\nu + 1}\,\rd y}
\]
is examined in \cite{banasiak.lamb.langer:2013}, and results are obtained
on the analyticity of associated fragmentation semigroups defined on the
weighted spaces $L^1(\mathbb{R}_+, (1 + x^m)\,\rd x)$;
see \cite[Theorems~2.1 and 2.3]{banasiak.lamb.langer:2013}.
\hfill$\lozenge$
\end{remark}

Motivated by the terms in \eqref{NA frag system}, we introduce, for each $t \in [0,T]$,
the formal expressions
\[
  \mathcal{A}(t): (f_n)_{n=1}^{\infty} \mapsto \bigl(-a_n(t)f_n\bigr)_{n=1}^{\infty}
\]
and
\[
  \mathcal{B}(t): (f_n)_{n=1}^{\infty} \mapsto
  \Biggl(\,\sum\limits_{j=n+1}^{\infty} a_j(t)b_{n,j}(t)f_j\Biggr)_{n=1}^{\infty}.
\]
Operator realisations, $A(t)$ and $B(t)$, of $\mathcal{A}(t)$ and $\mathcal{B}(t)$
respectively, are then defined in $\ell_w^1$ by
\begin{alignat}{2}
  A(t)f &= \mathcal{A}(t)f, \qquad &
  D\bigl(A(t)\bigr) &= \bigl\{f \in \ell_w^1: \mathcal{A}(t)f \in \ell_w^1\bigr\},
  \label{Adef}
  \\[1ex]
  B(t)f &= \mathcal{B}(t)f, \qquad &
  D\bigl(B(t)\bigr) &= D\bigl(A(t)\bigr).
  \label{Bdef}
\end{alignat}
That $B(t)$ is well defined on $D(A(t))$, for each $t \in [0,T]$,
can be seen as follows.
Assumption (A4) implies that, for $f=(f_n)_{n=1}^\infty\in(D(A(t)))_+$, we have
\begin{equation}\label{phiBphiA}
\begin{aligned}
  \phi_w\bigl(\mathcal{B}(t)f\bigr)
  &= \sum_{n=1}^\infty w_n\sum_{j=n+1}^\infty a_j(t)b_{n,j}(t)f_j
  = \sum_{j=2}^\infty\Biggl(\sum_{n=1}^{j-1}w_nb_{n,j}(t)\Biggr)a_j(t)f_j
  \\[1ex]
  &\le \sum_{j=2}^\infty \kappa w_j a_j(t)f_j
  \le -\kappa\sum_{j=1}^\infty w_j\bigl(-a_j(t)\bigr)f_j
  = -\kappa\phi_w\bigl(A(t)f\bigr);
\end{aligned}
\end{equation}
the change in the order of summation in the calculation above is justified since
each term is positive.  Now let $f=(f_n)_{n=1}^\infty\in D(A(t))$.
Then $|f|=(|f_n|)_{n=1}^\infty\in D(A(t))_+$, and we obtain from \eqref{phiBphiA} that
\begin{align*}
  \big\|\mathcal{B}(t)f\big\|_w
  &= \sum_{n=1}^\infty w_n\Bigg|\sum_{j=n+1}^\infty a_j(t)b_{n,j}(t)f_j\Bigg|
  \le \phi_w\bigl(\mathcal{B}(t)|f|\bigr)
  \\[1ex]
  &\le -\kappa\phi_w\bigl(A(t)|f|\bigr)
  = -\kappa\sum_{n=1}^\infty w_n\bigl(-a_n(t)\bigr)|f_n|
  = \kappa\big\|A(t)f\big\|_w < \infty,
\end{align*}
which yields $f\in D(B(t))$ and
\begin{equation}\label{eq NA B is bounded by A}
  \|B(t)f\|_w \le \kappa \|A(t)f\|_w \qquad \text{for all} \ f \in D\bigl(A(t)\bigr).
 \end{equation}
On setting $G(t) = A(t)+B(t)$,
we now write \eqref{NA frag system} as the non-autonomous ACP
\begin{equation}\label{eq NA frag ACP}
  u'(t) = G(t)u(t), \quad s < t \le T; \qquad u(s)=\mr{u}
\end{equation}
with $\mr{u}\in\ell_w^1$.

\medskip
The following results on the operators $G(t)$, $t\in[0,T]$, will be required
in the next section.

\begin{lemma}\label{analytic semigroup}
Let assumptions \textup{(A1)--(A4)} be satisfied.  For each $t\in[0,T]$,
\begin{myenuma}
\item
  the operator $G(t)$ is the generator of an analytic, substochastic $C_0$-semigroup,
  $(S_t(\tau))_{\tau\ge0}$, on $\ell^1_w$;
\item
  for $\lambda\in\CC$ with $\Re\lambda>0$, the resolvent operator $R(\lambda, G(t))$
  can be factorised as
  \begin{equation}\label{fact}
    R\bigl(\lambda,G(t)\bigr)
    = R\bigl(\lambda,A(t)\bigr)\Bigl[I-B(t)R\bigl(\lambda,A(t)\bigr)\Bigr]^{-1}
  \end{equation}
  where the factors on the right-hand side satisfy
  \begin{equation}\label{norm_est}
    \big\|R\bigl(\lambda,A(t)\bigr)\big\| \le \frac{1}{|\lambda|}\,, \qquad
    \bigg\|\Bigl[I-B(t)R\bigl(\lambda,A(t)\bigr)\Bigr]^{-1}\bigg\| \le \frac{1}{1-\kappa}\,.
  \end{equation}
\end{myenuma}
\end{lemma}

\begin{proof}
Part (a) is an immediate consequence of \cite[Theorem~5.2]{kerr2020discrete}.

For (b) let $\lambda\in\CC$ with $\Re\lambda>0$.
Since $a_n(t)\ge0$, we have $\lambda\in\rho(A(t))$ and, for $f\in\ell_w^1$,
\[
  \big\|R\bigl(\lambda,A(t)\bigr)f\big\|_w
  = \sum_{n=1}^\infty w_n\frac{1}{|\lambda+a_n(t)|}|f_n|
  \le \frac{1}{|\lambda|}\|f\|_w,
\]
which yields the first inequality in \eqref{norm_est}.
Observe that $D(B(t))=D(A(t))$ and hence
\begin{align}
  \lambda I-G(t) &= \lambda I-A(t)-B(t)
  = \Bigl(I-B(t)\bigl(\lambda I-A(t)\bigr)^{-1}\Bigr)
  \Bigl(\lambda I-A(t)\Bigr)
  \nonumber\\[1ex]
  &= \Bigl(I-B(t)R\bigl(\lambda,A(t)\bigr)\Bigr)\bigl(\lambda I-A(t)\bigr).
  \label{factorisation1}
\end{align}
For $f\in\ell_w^1$, we use \eqref{eq NA B is bounded by A} to obtain
\begin{align*}
  \big\|B(t)R\bigl(\lambda,A(t)\bigr)f\big\|_w
  &\le \kappa\big\|A(t)R\bigl(\lambda,A(t)\bigr)f\big\|
  \\[1ex]
  &= \kappa\sum_{n=1}^\infty w_n\frac{a_n(t)}{|\lambda+a_n(t)|}|f_n|
  \le \kappa\|f\|_w,
\end{align*}
which implies $\|B(t)R(\lambda,A(t))\|\le\kappa<1$.
Consequently, $I-B(t)R(\lambda,A(t))$ is invertible and
\[
  \bigg\|\Bigl[I-B(t)R\bigl(\lambda,A(t)\bigr)\Bigr]^{-1}\bigg\|
  = \Bigg\|\sum_{n=0}^\infty \Bigl[B(t)R\bigl(\lambda,A(t)\bigr)\Bigr]^n\Bigg\|
  \le \sum_{n=0}^\infty \kappa^n = \frac{1}{1-\kappa}\,,
\]
which proves the second inequality in \eqref{norm_est}.
Taking inverses on both sides of \eqref{factorisation1} we obtain \eqref{fact}.
\end{proof}

\section{Well-Posedness}
\label{sect:well-posedness}

We now establish sufficient conditions on the fragmentation coefficients and the
weight $w=(w_n)_{n=1}^\infty$ for the non-autonomous ACP \eqref{eq NA frag ACP}
to be well posed in the complex Banach lattice $\ell^1_{w}$.
In addition to requiring (A1)--(A4) to hold, we also assume that
constants $C_1\ge0$, $C_2\ge0$ and $\sigma\in(0,1]$ exist such that
\begin{myenum_A}
\setcounter{counter_c}{4}
\item
  for all $n\in\NN$ and $s,t,\tau\in[0,T]$,
  \[
    \frac{|a_n(t)-a_n(s)|}{1+a_n(\tau)} \le C_1|t-s|^\sigma;
  \]
\item
  for all $j\in\{2,3,\ldots\}$ and $s,t,\tau\in[0,T]$,
  \[
    \frac{1}{1+a_j(\tau)}
    \sum_{n=1}^{j-1} w_n\big|a_j(t)b_{n,j}(t)-a_j(s)b_{n,j}(s)\big|
    \le C_2w_j|t-s|^\sigma.
  \]
\end{myenum_A}

\begin{example}
Let $c_n,d_n\ge0$ for $n\in\NN$ and let $\varphi:[0,T]\to[K_1,\infty)$
be a function such that $|\varphi(t)-\varphi(s)|\le K_2|t-s|^\sigma$
with $K_1,K_2>0$.  Then
\[
  a_n(t) = c_n\varphi(t)+d_n
\]
satisfies (A5), which can be seen as follows: for $t,s,\tau\in[0,T]$ we have
\begin{align*}
  \frac{|a_n(t)-a_n(s)|}{1+a_n(\tau)}
  &= \frac{c_n|\varphi(t)-\varphi(s)|}{c_n\varphi(\tau)+d_n+1}
  \le \frac{K_2c_n|t-s|^\sigma}{c_nK_1+1}
  \\[1ex]
  &= \frac{K_2}{K_1}\cdot\frac{c_n}{c_n+\frac{1}{K_1}}|t-s|^\sigma
  \le \frac{K_2}{K_1}|t-s|^\sigma.
\end{align*}
\hfill$\lozenge$
\end{example}

\begin{remark}
If each $b_{n,j}$ is constant on $[0,T]$, say $b_{n,j}(t)\equiv b_{n,j}$, 
then (A6) follows from (A4) and (A5):
\begin{align*}
  & \frac{1}{1+a_j(\tau)}
  \sum_{n=1}^{j-1} w_n\big|a_j(t)b_{n,j}(t)-a_j(s)b_{n,j}(s)\big|
  = \sum_{n=1}^{j-1} w_n\frac{|a_j(t)-a_j(s)|}{1+a_j(\tau)}b_{n,j}
  \\[1ex]
  &\le C_1|t-s|^\sigma\sum_{n=1}^{j-1}w_n b_{n,j}
  \le C_1\kappa w_j|t-s|^\sigma.
\end{align*}
The case of constant $b_{n,j}$ corresponds to the situation when the outcome
of the fragmentation of an $n$-mer ($n\ge2$) does not depend on the time
at which it occurs.
For example, this arises in the Becker--D\"{o}ring model of a
coagulation--fragmentation process in which the break-up of an $n$-mer always results
in a monomer and an $(n-1)$-mer; see \cite[\S 2.2.1]{banasiak2019analytic}.
Note also that the coefficients $b_{n,j}$ in Example~\ref{example_powers} are constant on $[0,T]$.
\hfill$\lozenge$
\end{remark}

\begin{lemma}\label{constdomain}
Let assumptions \textup{(A1)--(A6)} hold.  Then,
\begin{myenuma}
\item
  $D(G(t)) = D(G(0)) \eqqcolon \mathcal{D}$ for all $t\in[0,T]$;
\item
  for all $s,t,\tau\in[0,T]$ we have
  \begin{equation}\label{est_Gt_Gs}
    \big\|\bigl(G(t)-G(s)\bigr)R\bigl(1,A(\tau)\bigr)\big\|
    \le (C_1+C_2)|t-s|^\sigma.
  \end{equation}
\end{myenuma}
\end{lemma}

\begin{proof}
(a)
Let $s,t\in[0,T]$ and assume that $f\in D(G(s))=D(A(s))$.
It follows from (A5) with $\tau=s$ that
\begin{align*}
  \sum_{n=1}^\infty w_na_n(t)|f_n|
  &\le \sum_{n=1}^\infty w_n\big|a_n(t)-a_n(s)\big|\,|f_n|
  + \sum_{n=1}^\infty w_n a_n(s)|f_n|
  \\[1ex]
  &\le \sum_{n=1}^\infty w_n C_1|t-s|^\sigma\bigl(1+a_n(s)\bigr)|f_n|
  + \sum_{n=1}^\infty w_n a_n(s)|f_n|
  \\[1ex]
  &= \bigl(C_1|t-s|^\sigma+1\bigr)\big\|A(s)f\|_w + C_1|t-s|^\sigma\|f\|_w
  < \infty,
\end{align*}
which implies that $f\in D(A(t))=D(G(t))$.
Since $s$ and $t$ were arbitrary, it follows that the domain of $G(t)$ is
independent of $t$.

(b)
Let $t,s,\tau\in[0,T]$ and $f\in\ell_w^1$.  Then
\begin{equation}\label{eqnew1}
\begin{aligned}
  \big\|\bigl(G(t)-G(s)\bigr)R\bigl(1,A(\tau)\bigr)f\big\|_w
  &\le \big\|\bigl(A(t)-A(s)\bigr)\bigl(I-A(\tau)\bigr)^{-1}f\big\|_w
  \\[1ex]
  &\quad + \big\|\bigl(B(t)-B(s)\bigr)\bigl(I-A(\tau)\bigr)^{-1}f\big\|_w.
\end{aligned}
\end{equation}
Let us estimate each term separately.  From (A5) we obtain
\begin{align}
  &\big\|\bigl(A(t)-A(s)\bigr)\bigl(I-A(\tau)\bigr)^{-1}f\big\|_w
  = \sum_{n=1}^\infty w_n\frac{|a_n(t)-a_n(s)|}{1+a_n(\tau)}|f_n|
  \nonumber\\[1ex]
  &\le \sum_{n=1}^\infty w_n C_1|t-s|^\sigma|f_n|
  = C_1|t-s|^\sigma\|f\|_w.
  \label{eqnew2}
\end{align}
For the second term on the right-hand side of \eqref{eqnew1}
we can use (A6) to deduce that
\begin{align}
  &\big\|\bigl(B(t)-B(s)\bigr)\bigl(I-A(\tau)\bigr)^{-1}f\big\|_w
  \nonumber\\[1ex]
  &= \sum_{n=1}^\infty w_n\Bigg|\sum_{j=n+1}^\infty
  \Bigl(a_j(t)b_{n,j}(t)-a_j(s)b_{n,j}(s)\Bigr)\frac{1}{1+a_j(\tau)}\,f_j\Bigg|
  \displaybreak[0]\nonumber\\[1ex]
  &\le \sum_{n=1}^\infty w_n\sum_{j=n+1}^\infty
  \frac{\big|a_j(t)b_{n,j}(t)-a_j(s)b_{n,j}(s)\big|}{1+a_j(\tau)}\,|f_j|
  \nonumber\\[1ex]
  &= \sum_{j=2}^\infty \sum_{n=1}^{j-1}w_n
  \frac{\big|a_j(t)b_{n,j}(t)-a_j(s)b_{n,j}(s)\big|}{1+a_j(\tau)}\,|f_j|
  \nonumber\\[1ex]
  &\le \sum_{j=2}^\infty C_2w_j|t-s|^\sigma|f_j|
  \le C_2|t-s|^\sigma\|f\|_w.
  \label{eqnew3}
\end{align}
Combining \eqref{eqnew1}, \eqref{eqnew2} and \eqref{eqnew3}
we arrive at \eqref{est_Gt_Gs}.
\end{proof}

To enable  Theorem \ref{thm pazy NA} to be applied, we rescale each
semigroup $(S_t(\tau))_{\tau\ge0}$ that is generated by $G(t)$ by setting
\[
  T_t(\tau) = e^{-\tau}S_t(\tau), \qquad \tau\ge0,\, t\in[0,T],
\]
and consider the associated non-autonomous ACPs
\begin{equation}\label{eq scaledNA frag ACP}
  v'(t) = H(t)v(t), \quad s < t \le T; \qquad v(s)=\mathring{v},
\end{equation}
for $s\in[0,T)$, where $H(t)=G(t)-I$ is the generator of $(T_t(\tau))_{\tau\ge0}$
for $t\in[0,T]$.

\begin{proposition}\label{prop solution of rescaled equation}
Let assumptions \textup{(A1)--(A6)} hold.
Then there exists an evolution family, $(V(t,s))_{0\le s \le t \le T}$ on $\ell_w^1$,
with the following properties:
\begin{myenuma}
\item
  $v(t) = V(t,s)\mr{v}$ is the unique classical solution in $\ell^1_w$
  of \eqref{eq scaledNA frag ACP} for any $\mathring{v} \in \ell_w^1$ and $s\in[0,T)$;
\item
  if $\mr v\ge0$, then $v(t) = V(t,s)\mr{v}\ge0$ for $t\in[s,T]$; \\
  if, in addition, $\mr v\ne0$, then $v(t)\ne0$ for $t\in[s,T]$.
\end{myenuma}
\end{proposition}

\begin{proof}
(a)
We show that the operators $H(t)$, $t\in[0,T]$, satisfy the assumptions (P1)--(P3)
of Theorem~\ref{thm pazy NA}.  As each $G(t)$ generates an analytic substochastic
$C_0$-semigroup on $\ell_w^1$, it follows immediately that each $H(t)$ is also the
generator of an analytic  $C_0$-semigroup on $\ell_w^1$.
Moreover, from Lemma~\ref{constdomain}\,(a), $D(H(t))=D(G(t))=\mathcal{D}$
for all $t\in[0,T]$.  Hence (P1) is satisfied.

For (P2), we use Lemma~\ref{analytic semigroup}\,(b) to obtain
\begin{equation}\label{eq2 resolvent bound}
  \big\|R\bigl(\lambda,H(t)\bigr)\big\| \le \frac{1}{(1-\kappa)|\lambda +1|} \qquad
  \text{for} \ \lambda\in\CC: \Re\lambda > -1.
\end{equation}
Let $\lambda = \alpha + i\beta$, where $\alpha \ge 0$ and $\beta\in\RR$.
On applying the arithmetic mean--quadratic mean inequality
we deduce that
\begin{align*}
  |\lambda+1| &= \sqrt{2}\cdot\sqrt{\frac{(\alpha+1)^2+\beta^2}{2}}
  \ge \sqrt{2}\cdot\frac{\alpha+1+|\beta|}{2}
  \\[1ex]
  &= \frac{\sqrt{\alpha^2+2\alpha|\beta|+|\beta|^2}+1}{\sqrt{2}}
  \ge \frac{\sqrt{\alpha^2+\beta^2}+1}{\sqrt{2}}
  = \frac{|\lambda|+1}{\sqrt{2}}\,.
\end{align*}
Together with \eqref{eq2 resolvent bound} we obtain
\[
  \big\|R\bigl(\lambda,H(t)\bigr)\big\|
  \le \frac{\sqrt{2}}{1-\kappa}\cdot\frac{1}{|\lambda|+1}
\]
for $\lambda\in\CC$ with $\Re\lambda\ge0$, which shows (P2).

Finally, for (P3), let $s,t,\tau\in[0,T]$.
From Lemmas~\ref{analytic semigroup}\,(b) and \ref{constdomain}\,(b) we can deduce that
\begin{align*}
  & \big\|\bigl(H(t)-H(s)\bigr)H(\tau)^{-1}\big\|
  = \big\|\bigl(G(t)-G(s)\bigr)R\bigl(1,G(\tau)\bigr)\big\|
  \\[1ex]
  &= \bigg\|\bigl(G(t)-G(s)\bigr)R\bigl(1,A(\tau)\bigr)
  \Bigl[I-B(\tau)R\bigl(1,A(\tau)\bigr)\Bigr]^{-1}\bigg\|
  \displaybreak[0]\\[1ex]
  &\le \Big\|\bigl(G(t)-G(s)\bigr)R(1,A(\tau)\bigr)\Big\|\,
  \bigg\|\Bigl[I-B(\tau)R\bigl(1,A(\tau)\bigr)\Bigr]^{-1}\bigg\|
  \\[1ex]
  &\le (C_1+C_2)|t-s|^\sigma\cdot\frac{1}{1-\kappa}.
\end{align*}
Therefore (P3) is also satisfied.
Now the assertion follows from Theorem~\ref{thm pazy NA}.

(b)
To establish that the unique classical solution, $v(t) = V(t,s)\mr{v}$,
is non-negative for all $t \in [s,T]$ whenever $\mr{v}\in(\ell_w^1)_+$,
we determine an infinite matrix representation of $V(t,s)$, $0 \le s \le t \le T$,
with respect to the natural Schauder basis $(e_n)_{n=1}^\infty$ for $\ell_w^1$
that is given by
\begin{equation}\label{basis}
  (e_n)_m = \begin{cases}
    1 \quad &\text{if} \ n=m,
    \\[0.5ex]
    0 \qquad &\text{otherwise}.
  \end{cases}
\end{equation}
Let $s\in[0,T)$ be fixed.
We begin by considering, for each fixed $n\in\NN$, the finite system of linear ODEs
\begin{equation}\label{finite system}
\begin{split}
  \frac{\partial}{\partial t}v_{m,n}(t,s) &= -\bigl(1+a_m(t)\bigr)v_{m,n}(t,s)
  + \sum_{j=m+1}^{n} a_j(t)b_{m,j}(t)v_{j,n}(t,s),
  \\
  &\hspace*{30ex} t \in(s,T], \;\; m = 1,2,\ldots,n;
  \\[1ex]
  v_{n,n}(s,s) &= 1; \qquad v_{m,n}(s,s) = 0, \quad m = 1,\ldots,n-1,
\end{split}
\end{equation}
where we set $\sum\limits_{j=n+1}^{n} a_j(t)b_{m,j}(t)v_{j,n}(t,s) = 0$.
It follows from assumptions (A5) and (A6) that all the coefficient functions,
$1+ a_m$ and $a_jb_{m,j}$, in \eqref{finite system} are continuous on $[0,T]$.
Standard ODE theory \cite[\S III.1]{hale2009ordinary} then establishes that,
for each $s$ and $n$, the system \eqref{finite system} has a unique solution.
If we now define
\begin{equation}\label{lower_triangle}
  v_{m,n}(t,s) \equiv 0, \qquad\text{for} \ m > n,
\end{equation}
then the resulting infinite sequence, $(v_{m,n}(t,s))_{m=1}^\infty$,
is a classical solution of the ACP \eqref{eq scaledNA frag ACP},
with $\mathring{v}=e_n$.
By uniqueness of classical solutions to \eqref{eq scaledNA frag ACP},
we can deduce that $V(t,s)e_n = (v_{m,n}(t,s))_{m=1}^\infty$.
Moreover, as $(e_n)_{n=1}^\infty$ is a Schauder basis for $\ell_w^1$,
and each operator $V(t,s)$ is linear and continuous on $\ell_w^1$, we obtain
\[
  V(t,s)f = \sum_{n=1}^{\infty} f_n V(t,s)e_n, \qquad f=(f_n)_{n=1}^\infty \in \ell_w^1,
\]
and hence, for each $m\in\NN$,
\[
  \bigl(V(t,s)f\bigr)_m = \Biggl(\sum_{n=1}^{\infty} f_n V(t,s)e_n\Biggr)_m
  = \sum\limits_{n=1}^{\infty}v_{m,n}(t,s)f_n,
\]
which can be interpreted as a matrix multiplication of the infinite matrix
$\mathbb{V}(t,s) = (v_{m,n}(t,s))_{m,n\in\NN}$ and $f=(f_n)_{n=1}^\infty$
written as a column vector.  It follows from \eqref{lower_triangle} that $\mathbb{V}(t,s)$
has the form
\begin{equation}\label{eq evolution family matrix}
  \mathbb{V}(t,s) =
  \begin{bmatrix}
    v_{1,1}(t,s) & v_{1,2}(t,s) & v_{1,3}(t,s) & \cdots\, \\[0.5ex]
    0 & v_{2,2}(t,s) & v_{2,3}(t,s) & \cdots\, \\[0.5ex]
    0 & 0 & v_{3,3}(t,s) & \cdots\, \\
    \vdots & \vdots & \vdots & \ddots\,
  \end{bmatrix}.
\end{equation}
Since $v_{m,n}$ are solutions of \eqref{finite system}, the entries of $\mathbb{V}(t,s)$
do not depend on the weight $w$.
In the following we show that all these entries are non-negative.
Let us start with the main diagonal.
For $m=n\in\NN$ the differential equation in \eqref{finite system} is
\[
  \frac{\partial}{\partial t}v_{n,n}(t,s) = -\bigl(1+a_n(t)\bigr)v_{n,n}(t,s), \quad t\in(s,T];
  \qquad v_{n,n}(s,s)=1,
\]
and therefore the terms in the leading diagonal of $\mathbb{V}(t,s)$ are given by
\begin{equation}\label{diagonal}
  v_{n,n}(t,s)
  = \exp\biggl(\!-\int_s^t \bigl(1+a_n(\tau)\bigr)\,\rd\tau\biggr)
  > 0.
\end{equation}
Next consider the case when $n>1$ and $m=n-1$.
Suppose that $v_{n-1,n}(t,s) < 0$
for $t$ in some maximal interval $(\varepsilon_{n-1},\hat{\varepsilon}_{n-1})$,
where $s \le \varepsilon_{n-1} < \hat{\varepsilon}_{n-1} \le T$.
From \eqref{finite system}, we have
\begin{equation}\label{eq NA k-1 equation}
\begin{aligned}
  \frac{\partial}{\partial t}v_{n-1,n}(t,s) &= -\bigl(1+a_{n-1}(t)\bigr)v_{n-1,n}(t,s)
  +a_n(t)b_{n-1,n}(t)v_{n,n}(t,s),
  \\
  &\hspace*{48ex} t \in (s,T].
\end{aligned}
\end{equation}
Since $v_{n,n}(t,s) \ge 0$ for $t \in [s,T]$, the right-hand side
of \eqref{eq NA k-1 equation} is positive on $(\varepsilon_{n-1},\hat{\varepsilon}_{n-1})$.
On the other hand, by continuity, $v_{n-1,n}(\varepsilon_{n-1},s) =0 $, and
therefore, by the Mean Value Theorem, there exists
$\varepsilon \in (\varepsilon_{n-1},\hat{\varepsilon}_{n-1})$
such that $\frac{\partial}{\partial t}v_{n-1}(t,s)\big|_{t=\varepsilon}<0$.
This is a contradiction, and so $v_{n-1,n}(t,s) \ge 0$
for all $t \in [s,T]$.

If $n>2$ and $m=n-2$, a similar argument shows that $v_{n-2,n}(t,s) \ge 0$ for $t \in [s,T]$,
and continuing in this way we obtain $v_{m,n}(t,s) \ge 0$ for all $t \in [s,T]$
and $m \le n$.
Since $v_{m,n}(t,s)\equiv 0$ for all $m>n$, it follows that $v_{m,n}(t,s)\ge 0$
for all $m,n\in\NN$ and $t\in[s,T]$.
From this it is immediate that $v(t)\ge0$ if $\mr v\ge0$.

To prove the last statement, let $\mr v\ne0$ so that $\mr v_m>0$
for some $m\in\NN$.  Then
\[
  v_m(t) = v_{m,m}(t,s)\mr v_m + \sum_{n=m+1}^\infty v_{m,n}(t,s)\mr v_n > 0
\]
for $t\in[s,T]$ since the first term is strictly positive by \eqref{diagonal}
and the infinite series is non-negative.
\end{proof}

Proposition~\ref{prop solution of rescaled equation} leads immediately to the
main result of the paper, namely the existence in $\ell_w^1$ of a unique, classical
solution of \eqref{eq NA frag ACP} and its positivity for positive initial conditions.
We also establish mass conservation under an additional assumption.
To this end, let us recall that the total mass is given by the
first moment, $M_1$, where
\begin{equation}\label{def_M1}
  M_1(f) \coloneqq \sum_{n=1}^\infty n f_n, \qquad f=(f_n)_{n=1}^\infty \in \ell_w^1.
\end{equation}

\begin{theorem}\label{thm NA solution of Frag ACP}
Let assumptions \textup{(A1)--(A6)} hold.  Then there exists
an evolution family $(U(t,s))_{0\le s \le t \le T}$ on $\ell_w^1$
such that the following statements are true.
\begin{myenuma}
\item
$u(t) = U(t,s)\mr{u}$ is the unique classical solution in $\ell^1_w$
of \eqref{eq NA frag ACP} for any $\mathring{u} \in \ell_w^1$.
\item
If $\mr{u}\ge0$, then the solution from \textup{(a)} satisfies $u(t)\ge0$ for $t\in[s,T]$; \\
if, in addition, $\mr u\ne0$, then $u(t)\ne0$ for $t\in[s,T]$.
\item
If
\begin{equation}\label{mass cons}
  a_1(t) = 0, \quad \sum_{n=1}^{j-1} nb_{n,j}(t) = j \qquad
  \text{for all} \ j=2,3,\ldots,\; t \in [0,T],
\end{equation}
then $M_1(u(t))=M_1(\mr u)$
for $t\in[s,T]$ and $\mr u\in(\ell_w^1)_+$.
\end{myenuma}
\end{theorem}

\begin{proof}
Let $(V(t,s))_{0\le s \le t \le T}$  be the evolution family on $\ell_w^1$ that
is associated with \eqref{eq scaledNA frag ACP}, and define
\[
  U(t,s) \coloneqq e^{t-s}V(t,s), \qquad 0 \le s \le t \le T.
\]
A routine argument shows that $(U(t,s))_{0\le s \le t \le T}$
is also an evolution family on $\ell_w^1$.  Further, for $\mr u\in\ell_w^1$, we have
\begin{align*}
  \frac{\partial}{\partial t}\Bigl(U(t,s)\mr u\Bigr)
  &= e^{t-s}V(t,s)\mr u + e^{t-s}\frac{\partial}{\partial t}\Bigl(V(t,s)\mr u\Bigr)
  \\[1ex]
  &= e^{t-s}\Bigl[V(t,s)\mr u+H(t)V(t,s)\mr u\Bigr]
  \\[1ex]
  &= e^{t-s}\bigl(I+H(t)\bigr)V(t,s)\mr u
  = G(t)U(t,s)\mr u.
\end{align*}
Hence $u(t) = U(t,s)\mr{u}$ is a classical solution of \eqref{eq NA frag ACP};
it is clearly non-negative on $[s,T]$ whenever $\mr{u} \in (\ell_w^1)_+$,
and it is non-zero when $\mr u\ne0$.
To establish that $u(t) = U(t,s)\mathring{u}$ is the unique classical solution
in $\ell_w^1$, we simply note that if another classical solution, say $\tilde{u}$, exists,
then \eqref{eq scaledNA frag ACP} has a second classical solution given
by $\tilde{v}(t) = e^{s-t}\tilde{u}(t)$, and this contradicts
Proposition~\ref{prop solution of rescaled equation}\,(a).

Finally, let us prove (c).
Since $M_1$ is a bounded linear functional on $\ell_w^1$ and the solution $u(t)$
and the coefficients $a_n(t)$, $b_{n,j}(t)$ are non-negative, we obtain
from \eqref{mass cons} that
\begin{align*}
  \frac{\rd}{\rd t}\Bigl[M_1\bigl(u(t)\bigr)\Bigr]
  &= M_1\bigl(u'(t)\bigr)
  = M_1\bigl(G(t)u(t)\bigr)
  \\[1ex]
  &= \sum_{n=1}^\infty n\Biggl(-a_n(t)u_n(t)+\sum_{j=n+1}^\infty a_j(t)b_{n,j}(t)u_j(t)\Biggr)
  \displaybreak[0]\\[1ex]
  &= -\sum_{n=1}^\infty na_n(t)u_n(t)
  + \sum_{j=2}^\infty\Biggl(\sum_{n=1}^{j-1} nb_{n,j}(t)\Biggr)a_j(t)u_j(t)
  \\[1ex]
  &= -\sum_{n=2}^\infty na_n(t)u_n(t) + \sum_{j=2}^\infty ja_j(t)u_j(t)
  = 0.
\end{align*}
From this we can deduce that $M_1(u(t))=M_1(\mr u)$ for $t\in[s,T]$.
\end{proof}

The classical solution that exists by Theorem~\ref{thm NA solution of Frag ACP}
satisfies the original infinite systems of equations \eqref{NA frag system}
since the $m$th component of $G(t)u(t)$ is the right-hand side of the $m$th equation
of \eqref{NA frag system}.
Note, however, that solutions of \eqref{NA frag system} are not unique in general;
see, for example, the discussion at the end of Section~4 in \cite{smith2012discrete}
for the autonomous case.

\section{Asymptotic Behaviour of Solutions}
\label{sect:asymptotics}

We now turn our attention to the long-time behaviour of classical solutions
to the non-autonomous ACP \eqref{eq NA frag ACP}, focussing on the mass-conserving case.
When \eqref{mass cons} holds, and the coefficients $a_n(t)$, $n\ge2$,
are strictly positive for all $t$, it is expected, from physical considerations,
that, if the unique solution $u(t) = U(t,s)\mr{u}$ exists for all $t \ge s$,
then $u(t)$ should converge to the monomeric state $M_1(\mr{u})e_1$ as $t\to\infty$.
There have been several related investigations into the asymptotic behaviour
of classical solutions to the autonomous ACP formulation of the constant-coefficient,
mass-conserving fragmentation system.
In particular, the expected convergence to $M_1(\mr{u})e_1$ is established in
the first moment space $X_{[1]}$ for constant-coefficient binary fragmentation
in \cite{carr1994asymptotic}, and for constant-coefficient multiple fragmentation
in \cite{banasiak2011irregular}.  The case of convergence in spaces $\ell^1_w$
for more general weights $w$ is discussed in \cite{banasiaklamb2012discrete},
where $w(x) = x^p,\,p>1$, and also in our recent paper \cite{kerr2020discrete}.
In both \cite{banasiaklamb2012discrete} and \cite{kerr2020discrete}, it is shown
that the convergence of solutions to the monomeric state is at an exponential rate,
which is given explicitly in \cite{kerr2020discrete}.
Our aim now is to adapt the arguments we used in \cite{kerr2020discrete}
to prove that, under suitable conditions on the time-dependent coefficients,
solutions to the mass-conserving non-autonomous fragmentation ACP, also converge
to a monomeric state at an explicitly defined exponential rate.
We begin with the following proposition.

\begin{proposition}\label{prop:decay}
Let assumptions \textup{(A1)--(A6)} hold, let $s\in[0,T)$ and set
\[
  \wt a_{s,T} \coloneqq \inf_{\tau\in[s,T]}\,\inf_{n\in\NN}a_n(\tau).
\]
Then
\begin{equation}\label{decayUts}
  \|U(t,s)\| \le \exp\bigl[-\wt a_{s,T}(1-\kappa)(t-s)\bigr],
  \qquad t\in[s,T],
\end{equation}
where $(U(t,s))_{0\le s \le t \le T}$ is the evolution family on $\ell_w^1$ whose
existence is established in Theorem~\ref{thm NA solution of Frag ACP}.
\end{proposition}

\begin{proof}
Let $\mr u\in(\ell_w^1)_+\setminus\{0\}$ be arbitrary, and let $u(t)=U(t,s)\mr u$
be the classical solution of \eqref{eq general NA ACP} from
Theorem~\ref{thm NA solution of Frag ACP}.
Further, let $\phi_w$ be defined as in \eqref{def_phi_w}.
Since $\phi_w$ is a bounded linear functional on $\ell_w^1$, we obtain from
\eqref{phiBphiA} that, for $\tau\in(s,T)$,
\begin{align}
  \frac{\rd}{\rd\tau}\phi_w\bigl(u(\tau)\bigr) &= \phi_w\bigl(u'(\tau)\bigr)
  = \phi_w\bigl(A(\tau)u(\tau)\bigr) + \phi_w\bigl(B(\tau)u(\tau)\bigr)
  \nonumber\\[1ex]
  &\le (1-\kappa)\phi_w\bigl(A(\tau)u(\tau)\bigr)
  = -(1-\kappa)\sum_{n=1}^\infty w_n a_n(\tau)u_n(\tau)
  \nonumber\\[1ex]
  &\le -(1-\kappa)\wt a_{s,T}\sum_{n=1}^\infty w_n u_n(\tau)
  = -(1-\kappa)\wt a_{s,T}\phi_w\bigl(u(\tau)\bigr).
  \label{diff_inequ}
\end{align}
By Theorem~\ref{thm NA solution of Frag ACP}\,(b), $u(\tau)\ge0$ and $u(\tau)\ne0$,
and hence $\phi_w(u(\tau))=\|u(\tau)\|>0$.
Dividing both sides of \eqref{diff_inequ} by $\phi_w(u(\tau))$ and
integrating over $\tau$ from $s$ to $t$ for $t\in(s,T]$ we deduce that
\[
  \phi_w\bigl(u(t)\bigr)
  \le \phi_w\bigl(u(s)\bigr)\exp\bigl[-\wt a_{s,T}(1-\kappa)(t-s)\bigr].
\]
Since $u(t)\ge0$ and $u(s)=\mr u\ge0$, this yields
\[
  \|U(t,s)\mr u\|_w = \|u(t)\|_w \le \|\mr u\|_w
  \exp\bigl[-\wt a_{s,T}(1-\kappa)(t-s)\bigr].
\]
It follows from the positivity of $U(t,s)$ and
\cite[Proposition~2.67]{banasiak2006perturbations} that
\[
  \|U(t,s)\| = \sup_{\mr u\ge0,\, \|\mr u\|_w\le1}\|U(t,s)\mr u\|_w
  \le \exp\bigl[-\wt a_{s,T}(1-\kappa)(t-s)\bigr],
\]
which is \eqref{decayUts}.
\end{proof}

\begin{remark}\label{rem:contraction}
Note that, in particular, $\|U(t,s)\|\le1$ for $s\le t$,
since $\wt a_{s,T}\ge0$ by assumption (A1).
\hfill$\lozenge$
\end{remark}

In Theorem~\ref{thm:decay_to_monomer} below we prove that, under certain assumptions,
the solution converges to a pure monomeric state as $t\to\infty$,
i.e.\ the state where only the first component is non-zero.
Let us therefore consider a decomposition of the space $\ell_w^1$ into
the span of $e_1$, where $e_1$ is defined in \eqref{basis}, and a complement.
Define the space
\[
  Y_w \coloneqq \biggl\{(f_n)_{n=2}^\infty: \sum_{n=2}^\infty w_n|f_n|<\infty\biggr\}
  \quad \text{with norm} \;\;
  \big\|(f_n)_{n=2}^\infty\big\|_{Y_w} \coloneqq \sum_{n=2}^\infty w_n|f_n|,
\]
and let
\begin{alignat*}{2}
  & J: Y_w \to \ell_w^1, \qquad & J\bigl[(f_n)_{n=2}^\infty\bigr] &= (0,f_2,f_3,\ldots),
  \\[1ex]
  & P: \ell_w^1 \to Y_w, \qquad & P\bigl[(f_n)_{n=1}^\infty\bigr] &= (f_2,f_3,f_4,\ldots),
\end{alignat*}
be the embedding of $Y_w$ into $\ell_w^1$ and the projection from $\ell_w^1$ onto $Y_w$
respectively.  Then $\ell_w^1=\spn\{e_1\}\oplus JY_w$.
Let us start with a little lemma.

\begin{lemma}\label{lem:decomp}
Let $w=(w_n)_{n=1}^\infty$ with $w_n\ge n$ for $n\in\NN$, and let $M_1$ be defined
as in \eqref{def_M1}.  Further, let $g\in\ell_w^1$ and
assume that $M_1(g)=0$.  Then
\begin{equation}\label{inequ_proj}
  \|g\|_w \le (w_1+1)\big\|Pg\big\|_{Y_w}.
\end{equation}
\end{lemma}

\begin{proof}
We can decompose $g=(g_n)_{n=1}^\infty$ as
\[
  g = g_1e_1 + J\wh g \qquad \text{with} \;\; \wh g=Pg=(g_n)_{n=2}^\infty.
\]
The assumption $0=M_1(g)=M_1(g_1e_1)+M_1(J\wh g)$ implies that
\[
  |g_1| = \big|M_1(g_1e_1)\big| = \big|{-M_1(J\wh g)}\big|
  \le \sum_{n=2}^\infty n|g_n| \le \sum_{n=2}^\infty w_n|g_n|
  = \big\|\wh g\big\|_{Y_w},
\]
which yields
\[
  \|g\|_w = w_1|g_1| + \big\|\wh g\big\|_{Y_w} \le (w_1+1)\big\|\wh g\big\|_{Y_w};
\]
this proves \eqref{inequ_proj}.
\end{proof}

The next proposition provides an explicit estimate for the distance of the solution
from a monomeric state on finite time intervals.
It is used in the proof of Theorem~\ref{thm:decay_to_monomer} below.

\begin{proposition}\label{prop:decay_to_monomer}
Let assumptions \textup{(A1)--(A6)} and relation \eqref{mass cons} hold.
Further, let $s\in[0,T)$ and set
\[
  \wh a_{s,T} \coloneqq \inf_{\tau\in[s,T]}\,\inf_{n\ge2}a_n(\tau).
\]
Let $\mr u\in(\ell_w^1)_+$ and let $u(t)=U(t,s)\mr u$ be the classical solution
of \eqref{eq general NA ACP} from Theorem~\ref{thm NA solution of Frag ACP}.
With $M_1$ defined as in \eqref{def_M1} we have
\begin{equation}\label{decay_monomer_prop}
  \big\|u(t)-M_1(\mr u)e_1\big\|_w
  \le (w_1+1)\|\mr u\|_w\exp\bigl[-\hat a_{s,T}(1-\kappa)(t-s)\bigr],
  \qquad t\in[s,T].
\end{equation}
\end{proposition}

\begin{proof}
Let us consider the matrix representation of $U(t,s)$ for $t\in[s,T]$, which is
obtained by multiplying the matrix in \eqref{eq evolution family matrix} with $e^{t-s}$,
and splitting  it according to the decomposition of the space into $\spn\{e_1\}$
and $JY_w$,
\[
  \mathbb U(t,s) =
  \left[\begin{array}{c|ccc}
    u_{1,1}(t,s) & u_{1,2}(t,s) & u_{1,3}(t,s) & \cdots\, \\[0.5ex]
    \hline
    0 & u_{2,2}(t,s) & u_{2,3}(t,s) & \cdots\rule{0ex}{3ex}\, \\[0.5ex]
    0 & 0 & u_{3,3}(t,s) & \cdots\, \\
    \vdots & \vdots & \vdots & \ddots\,
  \end{array}\right]
  \eqqcolon
  \begin{bmatrix}
    u_{1,1}(t,s) & U_{(12)}(t,s) \\[2ex]
    0 & U_{(22)}(t,s)
  \end{bmatrix}.
\]
Let $\mr u=(\mr u_n)_{n=1}^\infty\in(\ell_w^1)_+$ and fix $t\in(s,T]$.
The matrix representation $\mathbb{U}(t,s)$ yields
\begin{align*}
  g \coloneqq&\; U(t,s)\mr u - M_1(\mr u)e_1
  \\[1ex]
  =&\; \Bigl[u_{1,1}(t,s)\mr u_1+U_{(12)}(t,s)P\mr u-M_1(\mr u)\Bigr]e_1 + JU_{(22)}(t,s)P\mr u.
\end{align*}
By Theorem~\ref{thm NA solution of Frag ACP}\,(c) we have
\[
  M_1(g) = M_1\bigl(U(t,s)\mr u\bigr)-M_1(\mr u) = 0,
\]
which allows us to apply Lemma~\ref{lem:decomp} and obtain
\begin{equation}\label{est_g_proof}
  \|g\|_w \le (w_1+1)\|Pg\|_{Y_w}
  = (w_1+1)\big\|U_{(22)}(t,s)P\mr u\big\|_{Y_w}.
\end{equation}
In order to estimate the right-hand side, let us set
\[
  \wh w_n \coloneqq w_{n+1}, \qquad \wh a_n(t) \coloneqq a_{n+1}(t), \qquad
  \wh b_{n,j}(t) \coloneqq b_{n+1,j+1}(t)
\]
for $n,j\in\NN$.
It is easy to see that $\wh w_n$, $\wh a_n(t)$, $\wh b_{n,j}(t)$ satisfy
assumptions (A1)--(A6); for instance, (A6) can be checked as follows:
\begin{align*}
  & \frac{1}{1+\wh a_j(\tau)}\sum_{n=1}^{j-1} \wh w_n
  \big|\wh a_j(t)\wh b_{n,j}(t)-\wh a_j(s)\wh b_{n,j}(s)\big|
  \\[1ex]
  &= \frac{1}{1+a_{j+1}(\tau)}\sum_{n=2}^j w_n\big|a_j(t)b_{n,j}(t)-a_j(s)b_{n,j}(s)\big|
  \\[1ex]
  &\le \frac{1}{1+a_{j+1}(\tau)}\sum_{n=1}^j w_n\big|a_j(t)b_{n,j}(t)-a_j(s)b_{n,j}(s)\big|
  \\[1ex]
  &\le C_2 w_{j+1}|t-s|^\sigma
  = C_2\wh w_j|t-s|^\sigma.
\end{align*}
Since $u$ solves \eqref{eq general NA ACP},
the component $Pu(\cdot)=U_{(22)}(\,\cdot\,,s)P\mr u$
solves \eqref{eq general NA ACP} with $G(t)$ obtained by replacing
$w_n$, $a_n(t)$ and $b_{n,j}(t)$ by $\wh w_n$, $\wh a_n(t)$ and $\wh b_{n,j}(t)$
respectively.  Applying Proposition~\ref{prop:decay} to $U_{(22)}(t,s)$
and using \eqref{est_g_proof} we obtain
\begin{align*}
  \big\|u(t)-M_1(\mr u)e_1\big\|_w &= \|g\|_w
  \le (w_1+1)\|U_{(22)}(t,s)\|\,\|P\mr u\|_{Y_w}
  \\[1ex]
  &\le (w_1+1)\|\mr u\|_w\exp\bigl[-\wh a_{s,T}(1-\kappa)(t-s)\bigr],
\end{align*}
which proves \eqref{decay_monomer_prop}.
\end{proof}

In the next theorem, which is the main result of this section, we consider
solutions of the non-autonomous ACP
\begin{equation}\label{ACP_half-line}
  u'(t) = G(t)u(t), \quad  t\in(0,\infty); \qquad u(0)=\mr{u},
\end{equation}
where the operator $G(t)$ is defined for all $t\in[0,\infty)$.
We assume that $(w_n)_{n=1}^\infty$ and $\kappa$ are fixed so that
assumptions (A1)--(A6) hold for all $T\in(0,\infty)$.
It follows from Theorem~\ref{thm NA solution of Frag ACP}
that \eqref{ACP_half-line} has a unique classical solution in $\ell_w^1$
when $\mr u\in\ell_w^1$.

\begin{theorem}\label{thm:decay_to_monomer}
Let $a_n$ and $b_{n,j}$ be defined on $(0,\infty)$ for $n,j\in\NN$ and
let $w_n>0$, $n\in\NN$, and $\kappa\in(0,1)$ be such that assumptions \textup{(A1)--(A6)}
hold for every $T>0$ \textup{(}the constants $C_1,C_2,\sigma$ in \textup{(A5)},
\textup{(A6)} may depend on $T$\textup{)}.
Further assume that \eqref{mass cons} holds, let $\mr u\in(\ell_w^1)_+$,
and let $u$ be the unique classical solution of \eqref{ACP_half-line}.
\begin{myenuma}
\item
If
\[
  \wh a_{0,\infty} \coloneqq \inf_{t\in(0,\infty)}\,\inf_{n\ge2} a_n(t) > 0,
\]
then
\[
  \big\|u(t)-M_1(\mr u)e_1\big\|_w
  \le (w_1+1)\|\mr u\|_w\exp\bigl[-\wh a_{0,\infty}(1-\kappa)t\bigr],
  \qquad t\in[0,\infty).
\]
\item
If
\[
  \wh a \coloneqq \liminf_{t\to\infty}\,\inf_{n\ge2} a_n(t) > 0,
\]
then, for every $c<\wh a(1-\kappa)$ there exists $M>0$ such that
\begin{equation}\label{decay_to_monomer}
  \big\|u(t)-M_1(\mr u)e_1\big\|_w \le Me^{-ct}, \qquad t\in[0,\infty).
\end{equation}
\end{myenuma}
\end{theorem}

\begin{proof}
The assertion in (a) follows directly from Proposition~\ref{prop:decay_to_monomer}.

To prove (b), let $c<\wh a(1-\kappa)$.  There exists $s\in[0,\infty)$ such that
\[
  \frac{c}{1-\kappa} \le \wh a_{s,\infty}
  \coloneqq \inf_{t\in[s,\infty)}\,\inf_{n\ge2} a_n(t).
\]
It follows from Theorem~\ref{thm NA solution of Frag ACP}\,(c) and
Proposition~\ref{prop:decay_to_monomer} that, for $t\in[s,\infty)$,
\begin{align}
  \big\|u(t)-M_1(\mr u)e_1\big\|_w
  &= \big\|u(t)-M_1\bigl(u(s)\bigr)e_1\big\|_w
  \nonumber\\[1ex]
  &\le (w_1+1)\|u(s)\|_w \exp\bigl[-\wh a_{s,\infty}(1-\kappa)(t-s)\bigr]
  \nonumber\displaybreak[0]\\[1ex]
  &\le (w_1+1)\|\mr u\|_w\exp[-c(t-s)]
  \nonumber\\[1ex]
  &= (w_1+1)\|\mr u\|_w e^{cs}e^{-ct};
  \label{est1}
\end{align}
note that $\|u(s)\|_w\le\|\mr u\|_w$ by Remark~\ref{rem:contraction}.
For $t\in[0,s)$ we have
\begin{align*}
  \big\|u(t)-M_1(\mr u)e_1\big\|_w
  &\le (w_1+1)\|\mr u\|_w \exp\bigl[-\wh a_{0,s}(1-\kappa)t\bigr]
  \\[1ex]
  &\le (w_1+1)\|\mr u\|_w
  \le (w_1+1)\|\mr u\|_w e^{cs}e^{-ct},
\end{align*}
which, together with \eqref{est1} proves \eqref{decay_to_monomer}
with $M=(w_1+1)\|\mr u\|_w e^{cs}$.
\end{proof}

\section{Concluding Remarks}
\label{sect:concluding_remarks}

To summarise, in this paper we have used the theory of evolution families to
analyse the non-autonomous fragmentation system \eqref{NA frag system}.
By writing \eqref{NA frag system} as an ACP in an appropriately weighted $\ell^1$ space,
and exploiting results on the analyticity of semigroups associated with
autonomous fragmentation systems, obtained in our earlier paper \cite{kerr2020discrete},
we have proved the existence and uniqueness of classical solutions to the
non-autonomous problem, for time-dependent fragmentation coefficients that
satisfy the assumptions (A1)--(A6).
Properties of these solutions such as non-negativity and, under the additional
assumption \eqref{mass cons}, mass conservation have been established.
Moreover, results on the asymptotic behaviour of solutions have been obtained.

As mentioned in the Introduction, evolution families have also featured
in investigations into the non-autonomous continuous fragmentation equation,
which is given by
\begin{equation}\label{nonautcF}
\begin{split}
  \frac{\partial}{\partial t}u(x,t) &= -a(x,t)u(x,t)
  + \int_x^\infty a(y,t)b(x,y,t)u(y,t)\,\rd y,
  \\[1ex]
  &\hspace*{35ex} x\in(0,\infty),\; t\in(0,T],
  \\[1ex]
  u(x,0) &= \mathring{u}(x),
\end{split}
\end{equation}
where $u(x,t)$ represents the density of particles of size $x \in (0,\infty)$
at time $t$, and the coefficients $a(x,t)$ and $b(x,y,t)$ are interpreted in an
analogous manner to $a_n(t)$ and $b_{n,j}(t)$ in the discrete system \eqref{NA frag system}.
For the sake of comparison, we discuss briefly the key results that these
investigations have produced.

In \cite{mclaughlin1997NA}, a slightly different, but equivalent, formulation of the
initial-value problem \eqref{nonautcF} is posed as a non-autonomous ACP in the
space $L^1(\mathbb{R}_+, x\,\rd x)$ (denoted by $L_{1,-1}$ in \cite{mclaughlin1997NA}).
Only mass-conserving fragmentation is considered, and the fragmentation coefficients
are assumed to satisfy the following conditions:
\begin{myenum}
\item
for every $n>0$ there exists a function $C_n:[0,T]\to(0,\infty)$
such that
\begin{equation}\label{assump1}
  a(x,t) \le C_n(t), \qquad x\in(0,n],\; t\in[0,T];
\end{equation}
\item
there exists a function $G:(0,\infty)\times(0,\infty)\to(0,\infty)$ such that
\begin{equation}\label{assump2}
\begin{split}
  \big|a(y,t)b(x,y,t) - a(y,\tau)b(x,y,\tau)\big| \le |t-\tau|G(x,y), \hspace*{10ex} &
  \\[1ex]
  x,y\in(0,\infty),\; t,\tau\in[0,T], &
\end{split}
\end{equation}
where, for every $n>0$, $G$ is bounded on $(0,n]\times(0,n]$.
\end{myenum}
Under these assumptions, the existence of a strongly continuous
evolution family $(U(t,s))_{0 \le s \le t \le T}$, consisting of non-negative isometries
on $L^1(\mathbb{R}_+, x\,\rd x)$, is established.
Each  operator $U(t,s)$ is defined as the strong limit, as $n \to \infty$,
of operators $U_n(t,s)$, $n > 0$,  where, for each $n$, $(U_n(t,s))_{0\le s \le t \le T}$,
is a uniformly continuous evolution family that is associated with an appropriately
truncated version of \eqref{nonautcF},
where the truncation is with respect to $x$ to the interval $(0,n]$.
In the case of restricted initial data satisfying $\mathring{u}(x) \equiv 0$
on $[n,\infty)$, for some $n > 0$, it is shown that $u(t) = U(t,s)\mathring{u}$
is the unique classical solution of the non-autonomous ACP version of \eqref{nonautcF}.
However, there is no corresponding result for a general $\mr{u} \in L^1(\RR_+, x\,\rd x)$.
Instead, the function $u(t) = U(t,s)\mr{u}$ is interpreted as a `generalised' solution
of the non-autonomous ACP, and, provided  $C_n \in L^\infty([0,T])$,
where $C_n$ is the function in \eqref{assump1}, the associated scalar-valued
function $u(x,t) = [U(t,s)\mr{u}](x)$ is shown to be a solution of the
following integral version of \eqref{nonautcF}
\[
  u(x,t) = \mr{u}(x) - \int_{0}^t a(x,\tau)u(x,\tau)\,\rd\tau
  + \int_{0}^t\int_x^\infty a(y,\tau)b(x,y,\tau)u(y,\tau)\,\rd y\,\rd\tau.
\]
Some partial results on the uniqueness of the solution $u(t) = U(t,s)\mr{u}$,
for the case when $\mathring{u} \in L^1(\mathbb{R}_+, x\,\rd x)$ does not vanish
on $(n,\infty)$ for some $n > 0$, can be found in \cite{mclaughlin1997uniqueness},
where the notion of a weak solution is used.
In particular, it is shown that  $u(t) = U(t,s)\mr{u}$ is the unique, non-negative,
mass-conserving, weak solution of the non-autonomous ACP for any given non-negative
initial data $\mr{u} \in L^1(\RR_+,x\,\rd x)$,  provided that the
function $b$ is independent of time, and $a(x,t) = a_0(x)\alpha(t)$,
with $a_0(x) \le C_n$ on $(0,n]$, and $\alpha$ a Lipschitz continuous function on $[0,T]$.

More recently, evolution families, together with associated evolution semigroups,
have also been used in \cite{arlottibanasiak2010nonautonomous} to establish the
existence of a solution to the above integral version of \eqref{nonautcF},
still under the assumption that each fragmentation event conserves mass,
but with the milder restriction that the fragmentation rate $a$ only has to be
locally integrable with respect to time and locally bounded with respect to $x$.
As in \cite{mclaughlin1997NA}, the solution is given by $u(x,t) = [U(t,s)\mr{u}](x)$,
where $(U(t,s))_{0 \le s \le t \le T}$, is a strongly continuous evolution family of
non-negative contractive operators on $L^1(\RR_+,x\,\rd x)$.
Moreover, when $a$ is bounded on $[0,M] \times [0,T]$, for any $M,T\in(0,\infty)$,
each $U(t,s)$ is shown to be an isometry.

We believe that the approach we have used in this paper could also prove fruitful
if applied to appropriately posed ACP versions of \eqref{nonautcF}, and,
in particular, may lead to new results concerning the existence and uniqueness
of physically meaningful classical solutions.
A first step would clearly be that of identifying weighted spaces,
$L^1(\RR_+, w(x)\,\rd x)$, such that the semigroup associated with the autonomous,
continuous fragmentation equation is analytic when defined on $L^1(\RR_+, w(x)\,\rd x)$.
In connection with this, it is worth noting that sufficient conditions for the
fragmentation coefficients are stated in \cite[Section~5.1.7]{banasiak2019analytic}
which guarantee the analyticity of the continuous fragmentation semigroup for the
cases $w(x) = x^m$ and $w(x) = 1 + x^m$, where $m > 1$.

Finally, a natural extension of the work presented here is to incorporate coagulation
into the model, and then examine the full, non-linear, discrete
coagulation--fragmentation (C--F) system, in which both the coagulation and
the fragmentation coefficients are time-dependent.
Although an approach based on evolution families has been used
in \cite{mclaughlin1998existence}, for continuous C--F equations in which
the coagulation and fragmentation coefficients are both permitted to be time-dependent,
we are unaware of similar investigations into the discrete case.

\bigskip

\noindent
\textbf{Acknowledgements.} \\
L.~Kerr gratefully acknowledges the support of \textit{The Carnegie Trust
for the Universities of Scotland}.
Further, L.~Kerr is a cross-disciplinary post-doctoral fellow supported by 
funding from the \textit{University of Edinburgh} 
and \textit{Medical Research Council} (MC\_UU\_00009/2).


\noindent
Lyndsay Kerr: \\
MRC Institute of Genetics and Cancer, \\
University of Edinburgh, \\
United Kingdom. \\
\texttt{lyndsay.kerr@ed.ac.uk}

\vspace{2ex}

\noindent
Wilson Lamb, Matthias Langer: \\
Department of Mathematics and Statistics \\
University of Strathclyde \\
26 Richmond Street \\
Glasgow G1 1XH \\
United Kingdom \\
\texttt{w.lamb@strath.ac.uk, m.langer@strath.ac.uk}

\end{document}